\documentclass{tran-l}
\usepackage{amssymb,latexsym, amscd}
\usepackage[all]{xy}
\usepackage{graphicx}

 \newtheorem{theorem}{Theorem}[subsection]
 \newtheorem{cor}[theorem]{Corollary}
 \newtheorem{lemma}[theorem]{Lemma}
 \newtheorem{proposition}[theorem]{Proposition}
 \theoremstyle{definition}
 
 \theoremstyle{definition}
 
 \theoremstyle{remark}
 \newtheorem{rem}[theorem]{Remark}
 \numberwithin{equation}{subsection}

\usepackage{latexsym}
\usepackage{amssymb}
\newcommand{\ben}{\begin{equation}}
\newcommand{\een}{\end{equation}}

\newcommand{\integer}{\ensuremath{{\mathbb Z}}}

\newcommand{\real}{\ensuremath{{\mathbb R}}}
\newcommand{\complex}{\ensuremath{{\mathbb C}}}

\newcommand{\PP}{{\mathcal P}}

\newcommand{\LL}{\mathcal{L}}

\newcommand{\OO}{\mathcal{O}}

\newcommand{\Xx}{\mathsf{X}}

\newcommand{\Loop}{\mathsf{L}}

\newcommand{\gr}{\mathfrak}

\newcommand{\Rank}{\ensuremath{{\mathrm{Rank}}}}
\newcommand{\Map}{\ensuremath{{\mathrm{Map}}}}

\newcommand{\orb}{\ensuremath{{\mathrm{orb}}}}

\newcommand{\ev}{\ensuremath{{\mathrm{ev}}}}
\newcommand{\Symn}{\ensuremath{{\gr{S}_n}}}

\begin{document}

\title[The loop orbifold of the symmetric product]
{ The loop orbifold of the symmetric product}


\author{Ernesto Lupercio}
\address{Departamento de Matem\'{a}ticas, CINVESTAV,
     Apartado Postal 14-740
     07000 M\'{e}xico, D.F. M\'{E}XICO}
 \email{lupercio@math.cinvestav.mx}
 \author{Bernardo Uribe}    
\address{Departamento de Matem\'{a}ticas, Universidad de los Andes,
Carrera 1 N. 18A - 10, Bogot\'a, COLOMBIA}
\email{buribe@uniandes.edu.co}
\author{Miguel A. Xicotencatl}
\address{Departamento de Matem\'{a}ticas, CINVESTAV,
     Apartado Postal 14-740
     07000 M\'{e}xico, D.F. M\'{E}XICO}
\email{xico@math.cinvestav.mx}
\thanks{The second author was partially supported by the ``Fondo de apoyo
a investigadores jovenes" from Universidad de los Andes and part of his research was carried out at the MPIfM in Bonn. The first and third authors were partially supported by Conacyt.}

\date{June 15, 2006 and, in revised form, ------.}



\begin{abstract}
By using  the loop orbifold of the symmetric product, we give a formula for the Poincar\'e polynomial of the free loop space of the Borel construction of the symmetric product. We also show that the Chas-Sullivan product structure in the homology of the free loop space of the Borel construction of the symmetric product induces a ring structure in the homology of the inertia orbifold of the symmetric product. This ring structure is compared to the one in cohomology defined through the usual field theory formalism.
\end{abstract}

\maketitle

\tableofcontents

\section{Introduction}

The (naive) symmetric product of a space $X$ is often defined as the \emph{topological space} 
$$X^n/\Symn : = X\times\cdots\times X /\Symn.$$
We find that it is better to study instead the \emph{orbispace} 
$$[X^n/\Symn] : = [X\times\cdots\times X /\Symn].$$
Namely, the category whose objects are $n$-tuples $(x_1,\ldots,x_n)$ of points in $X$ and whose arrows are elements of the form $(x_1,\ldots,x_n; \sigma)$ where $\sigma \in \Symn$. The arrow $(x_1,\ldots,x_n; \sigma)$ has as its source $(x_1,\ldots,x_n)$, and as its target $(x_{\sigma(1)},\ldots,x_{\sigma(n)})$. This category is a groupoid for the inverse of $(x_1,\ldots,x_n; \sigma)$ is $(x_{\sigma(1)},\ldots,x_{\sigma(n)}; \sigma^{-1})$. For this reason we can think of $[X^n/\Symn]$ as an orbispace \cite{KontsevichTorus, Moerdijk2002}, and we call it the \emph{symmetric product of $X$}.

In this paper we study the basic properties of the topology of the loop orbispace of the symmetric product $[X^n/\Symn]$. By this we do not mean the free loopspace $\LL (X^n/\Symn)$ of the naive symmetric product, but rather the geometric realization of the loop orbispace $\Loop [X^n/\Symn]$, namely the free loopspace of the Borel construction $$Z_n : =\LL (X^n \times_\Symn E\Symn)= \Map(S^1, X^n \times_\Symn E\Symn).$$

Let us talk about the organization of this paper. In section 2 we collect some well-known facts about the symmetric product that will set the stage for what follows. In section 3 we prove the following formula for the  the generating function of the Poincar\'{e} polynomials of $Z_n$

\begin{theorem} Let $X$ be such that $H^i(\LL X;\real)$ is finitely generated. Let $\phi(Z_n,y)$ be the Poincar\'{e} polynomial of $Z_n$.
 Then
$$\sum_{n=0}^\infty \phi(Z_n,y) q^n =  \prod_{j >0} \frac{\prod_i (1 + q^j y^{2i+1})^{b^{2i+1}(\LL X)}}{\prod_i (1 - q^j y^{2i})^{b^{2i}(\LL X)}}$$ where $b^i(\LL X)$ is the $i$-th Betti number of $\LL X$.
\end{theorem}

Actually we prove a little bit more. We compute $H^*(\Loop [X^n/\Symn])$ with rational coefficients.

In section 4 we consider the case when $X=M$ is a smooth manifold and the orbifold is $\Xx=[M^n/\Symn]$. 

In \cite{LupercioUribeLoopGroupoid} we constructed a functor with image in infinite dimensional orbifolds $$\Loop \colon \mathbf{Orbifolds} \to \ \ \ S^1-\mathbf{Orbifolds}$$  so that when restricted to smooth manifolds it becomes the ordinary free loopspace functor $M\mapsto \LL M$. More interestingly the $S^1$ action on $\Loop \Xx$ has as a fixed suborbifold $I(\Xx)$ the inertia orbifold of $\Xx$. In \cite{deFxLuNevUri} we have argued that orbifold theories often localize to the inertia orbifold. 

Chas and Sullivan \cite{ChasSullivan} have defined an associative product on the homology of the loop space $H_*(\LL M)$. In \cite{LupercioUribeXicotencatl} we have generalized this construction from the category of manifolds to the category of orbifolds. In this paper we study this product on $H_*(\Loop [M^n/\Symn])$. Using this product and the localization principle mentioned above \cite{deFxLuNevUri} we define an associative product $(H_*(I[M^n/\Symn]),\bullet)$.

In section 5 we define a new product on the \emph{cohomology} of the inertia orbifold that we call the \emph{virtual intersection product} and we denote it by $\times$. To do this we use a criterion that Fantechi and G\"ottsche used to study the Chen-Ruan product on the cohomology of the inertia orbifold. Our definition is close to that of Chen and Ruan but we use the $d$ rather than the $\bar{\partial}$ operator to define our space of fields (the result is also a topological quantum field theory \cite{LupercioUribeQFT}).
 
We conclude the paper proving the following.

\begin{theorem} Under Poincar\'e duality we have the following ring isomorphism
$$  (H_*(I[M^n/\Symn]),\bullet) \cong (H^*(I[M^n/\Symn]),\times) $$
\end{theorem}

 {\it Acknowledgements.} The first and second authors would like to thank the hospitality of the Mathematical Sciences Research Institute in Berkeley. The second and third authors would like to thank the  Max Planck Institut in Bonn, and the University of Bonn, where part of this work was realized. We would like also to thank R. Cohen, N. Ganter, D. Gepner, A. Henriquez, I. Moerdijk, Y. Ruan, C. Teleman, A. Voronov, and C. Westerland for useful conversations.

\section{The symmetric product}
\subsection{Poincar\'e polynomials}
Let $X$ be a topological space, we will denote by $\phi(X,y)$ its Poincar\'e polynomial  $$\phi(X,y) = \sum_i b^i(X)y^i$$ where $b^i(X)$ is the $i$-th Betti number of $X$.

Macdonald \cite{Macdonald} proved the following formula,
$$\sum_{n=0}^\infty \phi(X^n/\gr{S}_n,y) q^n = \frac{\prod_i \left(1+qy^{2i+1}\right)^{b_{2i+1}(X)}}{\prod_i \left(1-qy^{2i}\right)^{b_{2i}(X)}}. $$

When we set the variable $y=-1$ we get the famous formula for the Euler characteristic of the symmetric product:
$$\sum_{n=0}^\infty \chi(X^n/\gr{S}_n) q^n = \left(1-q\right)^{-\chi(X)}. $$

The previous formul{\ae} are valid for topological spaces whose cohomology $H^i(X, \real)$ is finitely generated for each $i \geq 0$, and there is no restriction on the homological dimension of $X$.

\subsection{Equivariant (Orbifold) Euler characteristic}
There is a similar formula associated to the (equivariant) orbifold Euler characteristic $\chi_{\gr{S}_n}$ of the symmetric product, which is defined using the $\gr{S}_n$-equivariant $K$-theory of $X^n$ by the following expression,
$$\chi_{\gr{S}_n}(X^n) := \Rank\ K_{\gr{S}_n}^0(X^n) - \Rank\  K_{\gr{S}_n}^1(X^n)$$
and can also be calculated using generating functions by the following formula
\begin{eqnarray}
\sum_{n=0}^\infty \chi_{\gr{S}_n}(X^n) q^n = \prod_{j >0} \left(1 - q^j \right)^{-\chi(X)}. \label{equivariant Euler characteristic}
\end{eqnarray}

This last equation is obtained by using a formula due to Segal that allows to calculate the torsion free part of $K^*_G(Y)$ (where $G$ acts on $Y$ and $G$ is a finite group) by  localizing on the prime ideals of $R(G)$, the representation ring of $G$ \cite{Segal_K-theory}, namely $$K^*_G(Y) \otimes \complex \cong \bigoplus_{(g)} K^*(Y^g)^{C(g)}\otimes \complex$$ where $(g)$ runs over the conjugacy classes of elements in $G$, $Y^g$ are the fixed point loci of $g$ and $C(g)$ is the centralizer of $g$ in $G$.

For the symmetric group $\gr{S}_n$, its conjugacy classes are in one-to-one correspondence with partitions of $n$. Given  $ \tau \in \gr{S}_n$ we will write  $\sum_j jn_j =n$ to denote the  partition corresponding to its cojugacy class. Here $n_j$ stands for the number of cycles of size $j$ that appear in the $\tau$. The we have that the fixed point set $\left(X^n\right)^\tau$ is isomorphic to $X^{\sum_j n_j}$ and $C(\tau) \cong \prod_j \gr{S}_{n_j} \ltimes (\integer/j)^{n_j}$. As the cylic groups $\integer/ j$ act trivially in $K^*(X^{\sum_j n_j})$ the following decomposition holds
$$K^*_{\gr{S}_n}(X^n) \otimes \complex \cong \bigoplus_{(\tau)} K^*((X^n)^\tau)^{C(\tau)}\otimes \complex \cong \bigoplus_{\sum jn_j=n} \otimes_j K^*(X^{n_j})^{\gr{S}_{n_j}}\otimes \complex.$$

\subsection{Orbifold cohomology} For an orbifold $[Y/G]$ (viewed as a topological groupoid \cite{Moerdijk2002}) its {\it orbifold cohomology} is defined as the cohomology of the inertia orbifold $I[Y/G]$, i.e. $H_{orb}^*([Y/G]) := H^*(I[Y/G])$, where the inertia orbifold is defined as
$$I[Y/G] :=[ \left(\sqcup_{g \in G} Y^g\times \{g\} \right) / G]$$ where the action is given by
\begin{eqnarray*}
G \times \left(\sqcup_{g \in G} Y^g\times \{g\} \right) & \to & \left(\sqcup_{g \in G} Y^g\times \{g\} \right) \\
(h, (x,g)) & \mapsto & (xh, h^{-1}gh).
\end{eqnarray*}
 There is another presentation (Morita equivalent) of the inertia orbifold of $[Y/G]$ given by
 $$I[M/G] \cong \sqcup_{(g)} [Y^g/C(g)]$$
 where as before $(g)$ runs over the conjugacy classes, $Y^g$ is the fixed point loci and $C(g)$ is the centralizer.
Then we have $H_{orb}^*([Y/G]; \real) \cong \oplus_{(g)} H^*(Y^g;\real)^{C(g)}$, and by the chern character isomorphism
$K_G^*(Y) \otimes \complex \cong H_{orb}^*([Y/G];\complex)$.

We can define the Poincar\'e orbifold polynomial $\phi_{orb}([Y/G],y) = \sum b^i_{orb}([Y/G]) y^i$ where the orbifold Betti number $b^i_{orb}([Y/G])$ is the rank of $H_{orb}^i([Y/G]; \real)$.

For the symmetric product, viewed as an orbifold groupoid $[X^n/\gr{S}_n]$, we get that
\begin{eqnarray}
H_{\orb}^*([X^n/\gr{S}_n];\real) \cong \bigoplus_{\sum jn_j=n} \bigotimes_j H^*(X^{n_j};\real)^{\gr{S}_{n_j}} \label{orbifold cohomology of the symmetric product}
\end{eqnarray}
and calculating the orbifold Poincar\'e polynomial one gets
\begin{eqnarray} \label{equivariant Euler characteristic of symmetric product}
\sum_{n=0}^\infty \phi_{\orb}([X^n/\gr{S}_n],y) q^n & = & \sum_{n=0}^\infty q^n \left(\sum_{\sum jn_j=n } \prod_j \phi(X^{n_j}/\gr{S}_{n_j},y) \right) \\
& = & \sum_{n=0}^\infty  \left(\sum_{\sum jn_j=n } \prod_j \phi(X^{n_j}/\gr{S}_{n_j},y)(q^j)^{n_j} \right) \\
 & = & \prod_{j>0} \left( \sum_{n=0}^\infty \phi(X^n/\gr{S}_n,y) q^{jn} \right)\\
 &=& \prod_{j >0} \frac{\prod_i (1 + q^j y^{2i+1})^{b^{2i+1}(X)}}{\prod_i (1 - q^j y^{2i})^{b^{2i}(X)}}
\end{eqnarray}
that when $y=-1$, yields the formula \ref{equivariant Euler characteristic} for the equivariant Euler characteristic.

Again, for the previous formul{\ae} to be valid one only needs that the cohomology of $X$ is finitely generated at each $i$.

\begin{rem}
In algebraic geometry the orbifold cohomology is defined on the inertia orbifold but has a shift in grading, which is called {\it age} by M. Reid \cite{Reid} , {\it shifting number} by Chen-Ruan \cite{ChenRuan} and {\it fermionic shift} by physicists. Here we do not change the grading.
\end{rem}

\section{Loop orbifold of the symmetric product}

For an orbifold $[Y/G]$ the loop orbifold $\Loop [Y/G]$ has been defined in \cite{LupercioUribeLoopGroupoid, LupercioUribeXicotencatl} and for the case of a global quotient it has a very simple description: $\Loop [Y/G] = [\PP_GY/G]$ where $\PP_GY = \sqcup_{g \in G} \PP_gY \times \{g\}$ with $\PP_g Y = \{ f \colon [0,1] \to Y | f(0)g = f(1) \}$ and the $G$ action is given by
\begin{eqnarray*}
G \times \sqcup_{g\in G} \PP_gY \times\{g\} &\to&  \sqcup_{g\in G} \PP_gY \times\{g\} \\
(h, (f,g)) & \mapsto & (f \cdot h, h^{-1}gh)
\end{eqnarray*}
with $f\cdot h(t) := f(t)h$.
The loop orbifold has another presentation (Morita equivalent) given by
$$\Loop [Y/G] \cong \bigsqcup_{(g)} [\PP_gY/C(g)]$$
where $C(g)$ acts on $\PP_gY$ in the natural way.
It is a theorem proved in \cite{LupercioUribeXicotencatl} that  $B \Loop [Y/G] \simeq \LL B [Y/G]$, i.e. the geometrical realization of the loop orbifold is homotopically equivalent to the free loop space of the geometrical realization of the orbifold, which in terms of the Borel construction states: $$\bigsqcup_{(g)} \left(\PP_gY \times_{C(g)} EC(g) \right) \simeq \Map(S^1, Y\times_G EG).$$

For the case of the symmetric product, one gets $$\Loop [X^n/\gr{S}_n] \cong \bigsqcup_{(\tau)} [\PP_\tau X^n / C(\tau)].$$ But there is a better presentation of this orbifold, namely,
.

\begin{lemma}
The orbifold $[\PP_\tau X^n /C(\tau)]$ is isomorphic to $\prod_j [(\LL X)^{n_j} / \gr{S}_{n_j} \ltimes (\integer/j )^{n_j}]$ where the action of $\integer/ j $ is given by rotation by the angles $2\pi k/ j$ on $\LL X$, the free loop space of $X$.
\end{lemma}

\begin{proof}

When $(\tau)$ is represented by the product $\tau^1_1 \dots \tau^{n_1}_{1} \tau^1_2 \dots \tau^{n_2}_{2}\dots$ of disjoint cycles, with $\tau^i_j$ the $i$-th cycle of size $j$, and $\sum jn_j=n$, then $$\PP_\tau X^n \cong \prod_j \prod_{i=1}^{n_j} \PP_{\tau^i_j} X^j \cong \prod_j (\PP_{\sigma_j} X^j)^{n_j}$$
where $\sigma_j$ is the cycle $(1,2,\dots,j).$ Now, the space $\PP_{\sigma_j} X^j$ consists of $j$-tuples $f=(f_1, \dots, f_j)$ of paths $f_i \colon [0,1] \to X$ such that $f(0)\sigma_j = f(1)$, i.e. $f_i(0)=f_{\sigma_j(i)}(1)$, which imply that the paths $f_i$ could be concatenated into a loop $\tilde{f}$ which belongs to $\LL X$. The map $\PP_{\sigma_j} X^j \to \LL X$, $f \mapsto \tilde{f}$ is clearly a homeomorphism.

We have then,
$$[\PP_\tau X^n /C(\tau)] \cong \prod_j [(\PP_{\sigma_j} X^j)^{n_j} / \gr{S}_{n_j} \ltimes (\integer / j )^{n_j} ]
\cong \prod_j [(\LL X)^{n_j} / \gr{S}_{n_j} \ltimes (\integer / j )^{n_j}]$$
where the action of $\integer / j $ on an element $f=(f_1, \dots , f_j) \in \PP_{\sigma_j} X^j$ is generated by the action of $\sigma_j$, namely $f\cdot  \sigma_j= (f_j, f_1, \dots ,f_{j-1})$. As $f_j(0) = f_1(1)$, then the cyclic action rotates the loop $\tilde{f}$ by an angle of $2\pi/j$.
\end{proof}

As the action of $\integer/ j $ in $\LL X$ factors through the rotation action of the circle $S^1$ in $\LL X$, then

\begin{cor}
$$H^*(\Loop[X^n/\gr{S}_n] ;\real) \cong \bigoplus_{(\tau)} H^*(\PP_\tau X^n ;\real)^{C(\tau)} \cong \bigoplus_{\sum jn_j=n} \prod_j H^*((\LL X)^{n_j} ; \real)^{\gr{S}_{n_j}}$$
\end{cor}

At this point we can see some similarities between the loop orbifold of the symmetric product of $X$, and the inertia orbifold of the symmetric product of $\LL X$, namely that their rational cohomologies agree even though the orbifolds cannot be isomorphic

\begin{proposition}
The orbifolds $\Loop [X^n/\gr{S}_n]$ and $I [(\LL X)^n/\gr{S}_n]$ cannot be naturally isomorphic unless $n=1$, but their cohomologies with real coefficients agree.
\end{proposition}

\begin{proof}
By formula \ref{orbifold cohomology of the symmetric product} we have
$$H^*_{\orb}( [(\LL X)^n/\gr{S}_n] ;\real) \cong \bigoplus_{\sum jn_j=n} \prod_j H^*((\LL X)^{n_j} ; \real)^{\gr{S}_{n_j}}$$ which is isomorphic by the previous corollary to $H^*(\Loop[X^n/\gr{S}_n] ;\real)$.

But, the orbifolds $\Loop [X^n/\gr{S}_n]$ and $I [(\LL X)^n/\gr{S}_n]$ cannot be naturally isomorphic because the actions of the cyclic groups $\integer/ j $ are different. On the one hand, for $\Loop [X^n/\gr{S}_n]$, we just argued that the action of the cyclic groups are by rotation on $\LL X$ (coming from the action of $\sigma_j$ into $\PP_{\sigma_j} X^j$), and on the other, for $I [(\LL X)^n/\gr{S}_n]$, the action of the cyclic groups are trivial, because the copies of $\LL X$ come from the fixed point loci of the group action generated by the cycle $\sigma_j$ into $(\LL X)^j$. Therefore on the one hand one has the orbifold $[\LL X / (\integer/j)]$ with the rotation action, and in the other one has the orbifold $[\LL X / (\integer/j)]$ with the trivial action. These orbifolds cannot be naturally isomorphic.
In the case that $n=1$ both orbifolds are the same.

Let's see the case when $X=S^1$ and $n=2$. Then $\Loop [(S^1)^2/\gr{S}_2] = [(\LL S^1)^2 /\gr{S}_2] \sqcup [\LL S^1 / (\integer/2)]$ where the action of $\integer/ 2 $ in the second component is by rotation, and
 $ I[(\LL S^1)^2/\gr{S}_2] = [(\LL S^1)^2 /\gr{S}_2] \sqcup [\LL S^1 / \integer/2]$ where the action of $\integer / 2 $ is the trivial one. As $\LL S^1 \simeq \integer \times S^1$ it is easy to see that in the first case the geometrical realization of $[\LL S^1 / (\integer/2)]$ is homotopically equivalent to $(\integer \times S^1) \sqcup (\integer \times S^1 \times \real P^\infty)$ and in the second case is just $\integer \times S^1 \times \real P ^\infty$.
\end{proof}

Using the previous result and formula \ref{equivariant Euler characteristic of symmetric product}, one gets

\begin{cor}
Let $X$ be such that $H^i(\LL X;\real)$ is finitely generated. Then
$$\sum_{n=0}^\infty \phi(\Loop[X^n/\gr{S}_n],y) q^n =  \prod_{j >0} \frac{\prod_i (1 + q^j y^{2i+1})^{b^{2i+1}(\LL X)}}{\prod_i (1 - q^j y^{2i})^{b^{2i}(\LL X)}}$$ where $b_i(\LL X)$ is the $i$-th Betti number of $\LL X$.
And via the chern character map we get
$$K^*_{\gr{S}_n}((\LL X)^n) \otimes \complex \cong H^*(\Loop [X^n/\gr{S}_n];\complex).$$
\end{cor}

\begin{rem} \label{rem_exa_S^2-1}
The fact that the cohomologies of $I[\LL X^n /\gr{S}_n] $ and $\Loop [X^n/\gr{S}_n]$ agree is a feature of the symmetric product. In general, for any orbifold $[Y/G]$, the cohomologies of $I[\LL Y / G]$ and $\Loop [Y/G]$ do not have to agree. Take for example the $\integer/2 $ action on $S^2$ by rotating $\pi$ radians along the $z$-axis.
$I[\LL S^2/\integer / 2] = [\LL S^2/\integer / 2] \sqcup [\LL (S^2)^\xi / \integer/2]$ where $\xi$ generates the group $\integer/2$, and therefore $\LL (S^2)^\xi$ is the set of two points, the north and the south pole. Hence $H^*(I[\LL S^2/\integer / 2] ;\real) \cong H^*(\LL S^2; \real) \oplus \real^{\oplus 2}$.
On the other hand $\Loop [S^2/\integer/2] =  [\LL S^2/\integer / 2] \sqcup [\PP_\xi S^2/\integer/2]$ with cohomology $H^*(\Loop [S^2/\integer/2] ;\real) \cong H^*(\LL S^2 ; \real) \oplus H^*(\LL S^2 ;\real)$ (this is shown in the examples of \cite{LupercioUribeXicotencatl}).

\end{rem}

\section{Ring structure in the homology of the loop orbifold}
In \cite{LupercioUribeXicotencatl} we have showed that for orbifolds of the type $[M/G]$ with $M$ oriented, smooth and compact, and $G$ acting by orientation preserving diffeomorphisms, the homology of the loop orbifold  $H_*(\Loop[M/G])$ has the structure of a Batalyn-Vilkovisky algebra, i.e. a graded commutative algebra, with a degre $1$ operator $\Delta$, with $\Delta^2=0$,  and a Lie bracket that measures the discrepancy of $\Delta$ from being a derivation of the product.

 In this section we will study the ring structure of $H_*(\Loop[M^n/\gr{S}_n])$, and we will show that it induces a ring structure in the homology of   $I[M^n/\gr{S}_n]$ in such a way that its homology $H_*(I[M^n/\gr{S}_n])$  becomes a sub ring of $H_*(\Loop[M^n/\gr{S}_n])$.

So, let's start by showing the previous statement for $M$ itself

\begin{lemma}
The natural inclusion $i : M \to \LL M$ of constant loops and the evaluation at $0$, $ev : \LL M \to M$ induce ring maps in homology $i_* : H_*(M) \to H_*(\LL M)$  and $ev_* : H_*(\LL M) \to H_*(M)$ such that $ev_* \circ i_* = id$, in paticular  as $i_*$ is injective, $H_*(M)$ can be seen as a subring of $H_*(\LL M)$.
\end{lemma}
\begin{proof}
One just need to check that the following diagram is commutative
$$
\xymatrix{
\LL M \times_M \LL M  \ar[d]^{\ev_\infty} \ar[r] & \LL M \times \LL M \ar[d]^{ev \times ev}\\
M \ar@/^1pc/[u]^{i} \ar[r]^{diag} & M \times M. \ar@/^1pc/[u]^{i \times i}
}
$$
This induces the following diagram relating the Thom-Pontryagin construction of the top row with the bottom row (recall that the normal bundle of the diagonal inclusion is isomorphic to the tangent bundle, and the subindex $0$ means that we are taking everything ouside the zero section)
$$
\xymatrix{
\LL M \times \LL M \ar[d]^{ev \times ev} \ar[r] & (ev_\infty^* TM, (ev_\infty^* TM)_0) \ar[d]^{ev} \\
M\times M \ar@/^1pc/[u]^{i\times i} \ar[r] & (TM, TM_0) \ar@/^1pc/[u]^{i}
}
$$
that at the level of homology gives
$$
\xymatrix{
H_*(\LL M \times \LL M) \ar[d]^{ev_* \times ev_*}  \ar[r] & H_*(ev_\infty^* TM, (ev_\infty^* TM)_0) \ar[d]^{ev_*} \ar[rr]^\cong &  &H_{*-d}(\LL M) \ar[d]^{ev_*} \\
H_*(M \times M) \ar@/^1pc/[u]^{i_*\times i_*} \ar[r] & H_*(TM, TM_0) \ar[rr]^= & & H_{*-d}(M) \ar@/^1pc/[u]^{i_*}
}
$$ where $d=\mbox{dim}(M)$.
Then one has that $i_*$ and $ev_*$ are ring homomorphism, and as $ev \circ i = id$ then $i_*$ is injective
\end{proof}

For the case of the loop orbifold of the symmetric product, let's recall from \cite{LupercioUribeXicotencatl} how is the ring structure defined.
As the following diagram is a pull-back square
$$\xymatrix{
\PP_\tau M^n {}_1\times_0  \PP_\sigma M^n \ar[r] \ar[d]^{ev_{\infty}} & \PP_\tau M^n \times  \PP_\sigma M^n \ar[d]^{ev_1 \times ev_0} \\
M^n \ar[r] & M^n \times M^n,
}$$
one can do the Thom-Pontryagin construction, defining a homomorphism
$$H_*(\PP_\tau M^n \times \PP_\sigma M^n) \to H_{*-nd}(\PP_{\tau\sigma} M^n)$$
where the map $H_*(\PP_\tau M^n {}_1 \times_0 \PP_\sigma M^n) \to H_{*}(\PP_{\tau\sigma} M^n)$ is induced by the natural concatenation of paths $\PP_\tau M^n {}_1 \times_0 \PP_\sigma M^n \to \PP_{\tau\sigma} M^n$.

Then we have a product
\begin{eqnarray*}
H_p(\PP_\tau M^n)  \times H_q( \PP_\sigma M^n) & \to & H_{p+q-nd}(\PP_{\tau\sigma} M^n) \\
(\alpha, \beta) & \mapsto & \alpha \cdot \beta
\end{eqnarray*}
that is graded (shifted by $-nd$)  associative, and thus defines a product in $$ \bigoplus_\tau H_*(\PP_\tau M^n) \times \{\tau\}.$$
By taking the $\gr{S}_n$ invariant part
$$  \left( \bigoplus_\tau H_*(\PP_\tau M^n) \times \{\tau\} \right)^{\gr{S}_n} \cong H_*(\Loop [ M^n / \gr{S}_n])$$
we have defined thus a ring structure in the homology of the loop orbifold of the symmetric product.

Now let's study what is the behavior of the evaluation and inclusion of constant maps. So
consider the following commutative diagram
$$
\xymatrix{
\PP_\tau M^n \ar[r]^{ev} & M^n \\
(M^n)^\tau \ar[ru]_{f^\tau} \ar[u]^{i^\tau}
}
$$
where $f^\tau$ is the lnclusion of fixed point set, $i^\tau$ is the inclusion of constant loops, and $ev$ is the evaluation at $0$, we have the following
\begin{lemma}
The image in homology of $ev_*$ is equal to the image in homology of $f^\tau_*$.
\end{lemma}
\begin{proof}
Restricting the previous diagram to one of the cycles $\sigma$ of size $l$ that defines $\tau$, the diagram becomes
$$
\xymatrix{
\PP_\sigma M^l = \LL M \ar[rr]^{ev} && M^l \\
(M^l)^\sigma = M \ar[rru]_{f^\sigma} \ar[u]^{i^\sigma}
}
$$
where $f^\sigma$ becomes the diagonal inclusion $M \to M^l$ and the evaluation map $ev$ takes a loop $\alpha : S^1 \to M$ and maps it to $ev(\alpha) = (\alpha(0), \alpha(\frac{2\pi}{l}), \dots, \alpha(\frac{2(l-1)\pi}{l}))$. Defining the homotopy $ev^t(\alpha) = (\alpha(0), \alpha(\frac{2\pi t}{l}), \dots, \alpha(\frac{2(l-1)\pi t}{l}))$ one sees that $ev^1=ev $ and $ev^0$ are homotopic, and as $ev^0(\alpha)= f^\sigma(\alpha(0))$, the lemma follows.
\end{proof}

As the inclusion maps $f^\tau$ induce injective homomorphisms $f^\tau_* : H_*((M^n)^\tau) \to H_*(M^n)$, we define the groups $H^\tau_*(M^n):= image(f^\tau_*) \subset H_*(M^n)$ that with the use of the previous lemma, we get
$$
\xymatrix{
H_*(\PP_\tau M^n) \ar[r]^{ev_*} & H^\tau_*( M^n) \\
H_*((M^n)^\tau) \ar[ru]^\cong_{f^\tau_*} \ar[u]^{i^\tau_*}
}
$$

So we can define a ring structure in $\bigoplus_\tau H^\tau_*(M^n) \times \{\tau\}$ in the following way
\begin{eqnarray*}
\bullet : (H^\tau_*(M^n)\times\{\tau\} )\times (H^\sigma_*(M^n)\times \{ \sigma\})& \to & (H^{\tau\sigma}_{*-nd} (M^n) \times \{\tau\sigma\}) \\
((\alpha,\tau),( \beta, \sigma)) & \mapsto & (\alpha \bullet \beta, \tau\sigma)\end{eqnarray*}
where $$\alpha \bullet \beta = ev_*\left( \left( i^\tau_* \circ (f^\tau_*)^{-1} \alpha\right) \cdot  \left( i^\sigma_* \circ (f^\sigma_*)^{-1} \beta \right) \right)$$
and $\cdot$ is the product structure in the loop orbifold. Using the isomorphisms $f^\tau_*$ we also have a ring structure in $\bigoplus_\tau H_*((M^n)^\tau) \times \{\tau\}$ that we will also denote by $\bullet$.

Then we have the compatibility of all the products
$$\xymatrix{
H_*((M^n)^\tau) \times H_*((M^n)^\sigma) \ar@/^1pc/[rr]^{\cong} \ar[r]_{i^\tau_* \times i^\sigma_*} \ar[d]^\bullet &
H_*(\PP_\tau M^n) \times H_*(\PP_\sigma M^n) \ar[r]_{ev_* \times ev_*} \ar[d]^\cdot &
H^{\tau}_*(M^n) \times H^\sigma_*(M^n) \ar[d]^\bullet \\
H_*((M^n)^{\tau\sigma})  \ar[r]^{i^{\tau \sigma}_* }  \ar@/_1pc/[rr]_{\cong}  &
H_*(\PP_{\tau \sigma} M^n)  \ar[r]^{ev_*}  &
H^{\tau\sigma}_*(M^n)
}$$
and by taking $\gr{S}_n$ invariants we know  that
$$H_*(I[M^n/\gr{S}_n]) \cong \left( \bigoplus_\tau H^\tau_*(M^n) \times \{\tau\} \right)^{\gr{S}_n},  $$
so we can conclude
\begin{proposition}
The homology of the inertia orbifold $(H_*(I[M^n/\gr{S}_n]), \bullet)$ becomes an associative graded (with grading shifted by $-nd$) ring. Moreover, the inclusion of constant loops $i: I[M^n/\gr{S}_n] \to \Loop [M^n/\gr{S}_n]$  and the evaluation maps induce ring homomorphisms that makes the following diagram commute
$$\xymatrix{
& H_* (\Loop [M^n/\gr{S}_n]) \ar[dr]^{ev_*}  & \\
H_*(I[M^n/\gr{S}_n] )  \ar[ur]^{i_*} \ar[rr]_\cong && \left(\bigoplus_\tau H^\tau_*(M^n) \times \{\tau\} \right)^{\gr{S}_n}
}.$$
\end{proposition}

\begin{rem}
The inclusion of the inertia orbifold into the loop orbifold, in general does not induce an injective homomorphism in homology. Take the example of remark \ref{rem_exa_S^2-1}, namely the action of $\integer /2$ in $S^2$ by rotation along the $z$-axis. If the generator of $\integer /2$ is $\xi$, then the fixed point set $(S^2)^\xi$ consist of two points, the north and the south pole. The inclusion of the inertia orbifold into the loop orbifold is then $(S^2)^\xi \to \PP_\xi S^2$, where $\PP_\xi S^2 =\{ f :[0,1] \to S^2 | f(0)\xi = f(1) \}$. It is clear that $\PP_\xi S^2 \simeq \LL S^2$ which is connected, then the homomorphism $H_*((S^2)^\xi) \to H_*(\PP_\xi S^2)$ is not injective.
\end{rem}

\begin{rem} We have seen how to define a ring structure in the homology of $I[M^n/\gr{S}_n]$ using the structure of the homology of the loop orbifold. It is easy to see that the homology product we have defined boils down to intersection of cycles in $M^n$. Namely, for cycles  in $(M^n)^\tau$ and $(M^n)^\sigma$ (say $\alpha \in H^\tau_*(M^n)$ and $ \beta\in H^\sigma_*(M^n)$), their transversal intersection in $M^n$ is a cycle in $(M^n)^{\langle \tau, \sigma \rangle}$ ($\alpha \cap \beta \in H^{\tau,\sigma}_{*-nd}(M^n)$), and therefore could be pushforwarded to a cycle in $(M^n)^{\tau\sigma}$ ($\alpha \cap \beta \in H_{*-nd}^{\tau\sigma}(M^n)$). The associativity follows directly from the fact that transversal intersection is associative in homology.
\end{rem}

\section{The virtual intersection product of an orbifold.}

We would like to compare the product structure that we have defined in the previous section on the inertia orbifold to other products that exist on the same space, in particular the Chen-Ruan product \cite{ChenRuan}. For that purpose we are going to summarize a criterion of Fantechi and G\"otsche \cite{FantechiGotsche} on how to define a product in the cohomology of the inertia orbifold.

Consider the complex orbifold $[Y/G]$ where $Y$ is a complex manifold and $G$ acts holomorphically. Define the groups
$$H^*(Y,G) := \bigoplus_{g \in G} H^*(Y^g) \times \{g\}$$
where $Y^g$ is the fixed point set of the element $g$. The group $G$ acts in the natural way on the cohomologies and by conjugation on the labels. Denote by $Y^{g,h} = Y^g \cap Y^h$ and let's suppose we have $G$ invariant cohomology classes $c(g,h) \in H^*(Y^{g,h})$; i.e. such that $v^* c(k^{-1}gk,k^{-1}hk) = c(g,h)$ where $v : Y^{k^{-1}gk, k^{-1}hk} \to Y^{g,h}$ takes $x$ to $v(x):=xk$.
Now define the map
\begin{eqnarray*}
\times : H^*(Y^g) \times H^*(Y^h) & \to & H^*(Y^{gh}) \\
(\alpha, \beta) & \mapsto & i_* \left( \alpha|_{Y^{g.h}} \cdot \beta|_{Y^{g,h}} \cdot c(g,h) \right)
\end{eqnarray*}
where $i : Y^{g,h} \to Y^{gh}$ is the natural inclusion. By lemma 1.17 of \cite{FantechiGotsche}  a sufficient condition for the map $\times$ to define an associative product on $H^*(Y,G)$ is that, for every ordered triple of elements $(g,h,k) \in G$ the following relation holds in the cohomology of $W=Y^{g} \cap Y^{h} \cap Y^{k}$:
$$c(g,h)|_W \cdot c(gh,k)|_W \cdot e(Y^{gh}, Y^{g,h}, Y^{gh, k}) = c(g, hk)|_W \cdot c(h,k)|_W \cdot e(Y^{hk}, Y^{g,hk},Y^{h,k})$$
where $e(S,S_1,S_2)$ stands for the Euler class of the excess intersection bundle $E(S,S_1,S_2)$ when $S_1$ and $S_2$ are closed submanifolds of $S$. This bundle measures the failure of $S_1$ and $S_2$ to interesect transversally in $S$ having the property
$$j_2^* j_{1*} \alpha = i_{2*}(e(S,S_1,S_2) i_1^*(\alpha))$$
where $j_i : S_i \to S$ and $i_j : U \to S_j$ are the inclusions and $U=S_1 \cap S_2$.
 In the Grothendieck group of vector bundles the excess bundle is
$$E(S,S_1,S_2) = T_S|_U + T_U - T_{S_1}|_U - T_{S_2}|_U.$$

In particular, we have
\begin{lemma}
If $c(g,h)=e(Y,Y^g,Y^h)$ then $\times$ defines an associative product on $H^*(Y,G)$.
\end{lemma}
\begin{proof}
As $e(E + F) = e(E) e(F)$ we just need to check the equality in the Grothendieck ring of vector bundles over $W$.
The left hand side is
$$E(Y,Y^g,Y^h)|_W + E(Y, Y^{gh},Y^k)|_W + E(Y^{gh}, Y^{g,h} Y^{gh,k}) =  \ \ \ \ \ \ \ \ \ \ \ \ \ \ \  $$
$$T_Y + T_{Y^{g,h}} - T_{Y^g} - T_{Y^h} + T_{Y} +T_{Y^{gh,k}} -T_{Y^{gh}} -T_{Y^{k}} + T_{Y^{gh}} + T_{Y^{g,h,k}} - T_{Y^{g,h}} -T_{Y^{gh,k}}$$
(all the bundles are restricted to $W$) and after a reordering one can see that is equal to
$$E(Y,Y^g,Y^{hk})|_W + E(Y, Y^{h},Y^k)|_W + E(Y^{hk}, Y^{g,hk}, Y^{h,k})$$
\end{proof}

We will call the product $\times$  {\it virtual intersection product}. It is different  from the Chen-Ruan orbifold  product \cite{ChenRuan, FantechiGotsche, Uribe} because that one intersects holomorphically cycles in triples and therefore there is less room to do perturbation theory (in the virtual intersection product one has the operator $d$ on the moduli space, and in the Chen-Ruan product one has the operator $\bar{\partial}$). One can see this fact clearly, because the degree of the classes $c_{CR}(g,h)$ for the Chen-Ruan product is smaller or equal than the degree of the classes $e(Y, Y^g, Y^h)$ defined above.

In the case of the symmetric product
$$\mbox{deg} (e(M^n, (M^n)^\tau, (M^n)^\sigma)) = d[n + \OO(\langle \tau, \sigma \rangle) - \OO(\langle \tau \rangle) - \OO(\langle \sigma \rangle)]$$
where $d=dim_\real(M)$ and $\OO(\Gamma)$ is the number of orbits of the action of $\Gamma \subset \gr{S}_n$ into $\{1,2,\dots,n\}$, and
$$ \mbox{deg} ( c_{CR}(\tau, \sigma)) = \frac{d}{2}[n + 2\OO(\langle \tau, \sigma \rangle) - \OO(\langle \tau \rangle) - \OO(\langle \sigma \rangle) - \OO(\langle \tau\sigma \rangle) ],$$ see \cite{Uribe}.
As $\langle \tau\sigma \rangle$ is a subgroup of $\langle \tau,\sigma \rangle$, then $\OO(\langle \tau, \sigma \rangle) \leq \OO(\langle \tau \sigma \rangle) $, and therefore $$\mbox{deg} (e(M^n, (M^n)^\tau, (M^n)^\sigma)) \geq 2 \ \mbox{deg} ( c(\tau, \sigma)).$$

In the symmetric product, it is easy to see that the product $\times$ we have defined in the cohomology of the inertia orbifold is just the Poincar\'e dual of the product $\bullet$ in homology we defined previously. Using the isomorphisms $f_*^\tau : H_*((M^n)^\tau) \cong H_*^\tau(M^n)$ we have the following commutative diagram:
$$
\xymatrix{
H_p^\tau(M^n) \times H_q^\sigma(M^n) \ar@{<->}[rr]^{PD}  \ar[dd]^\cap  \ar@<-5ex>@/_2pc/[ddd]_{\bullet}
&& H^{d\OO(\langle \tau \rangle ) - p}((M^n)^\tau) \times H^{d\OO(\langle \sigma \rangle) - q}((M^n)^\sigma)  \ar[d]_{|_{(M^n)^{\langle \tau,\sigma\rangle }}}  \ar@<+12ex>@/^2pc/[ddd]^{\times}\\
  && H^{d\OO(\langle \tau \rangle ) + d\OO(\langle \sigma \rangle ) -p - q} ((M^n)^{\tau,\sigma}) \ar[d]_{\cup e} \\
  H_{p + q -nd}^{\tau,\sigma} (M^n) \ar@{>->}[d]^{inclusion}  \ar@{<->}[rr]^{PD} && H^{dn +d\OO(\langle \tau, \sigma \rangle) -p - q} ((M^n)^{\tau,\sigma}) \ar[d]_{pushforward}\\
H_{p + q -nd}^{\tau\sigma} (M^n) \ar@{<->}[rr]^{PD} && H^{dn +d\OO(\langle \tau \sigma \rangle ) -p - q} ((M^n)^{\tau,\sigma})
}
$$
where $PD$ denotes Poincar\'e duality.  Therefore we can conclude,
\begin{proposition}
The rings
$$\left( (\oplus_\tau H_*((M^n)^\tau) \times \{\tau\} ) , \bullet \right) \ \ \ \ \ \mbox{and} \ \ \ \ \ \
\left( (\oplus_\tau H^{d\OO(\tau)-*}((M^n)^\tau) \times \{\tau\} ) , \times \right)$$
are isomorphic under the Poincar\'e duality map. Therefore, taking $\gr{S}_n$ invariants we have then
$$(H_*(I[M^n/\gr{S}_n]) , \bullet)\cong (H^*(I[M^n/ \gr{S}_n]), \times).$$
\end{proposition}

Here it may be worthwhile to mention that the same theorems are valid if we use $K$-theory rather than cohomology \cite{JarvisKaufmannKimura}. The proofs are the same.

\bibliographystyle{amsplain}
\bibliography{Loop_orb_sym_pro}

\end{document}